\documentclass[11pt]{article}
\usepackage{epic,latexsym,amssymb}
\usepackage{color}
\usepackage{tikz}

\usepackage[format=hang,justification=centering]{caption}

\usetikzlibrary{shapes.geometric}
\tikzstyle{snode}=[circle,draw=black,fill=white,thick, inner sep=0pt ,minimum size=1.2mm]
\tikzstyle{gamma2}=[rectangle ,draw=black,inner sep=0pt ,minimum size=3.2mm]
\tikzstyle{bnode}=[circle ,draw=black,fill=black,thick, inner sep=0pt ,minimum size=1.2mm]
\tikzstyle{alpha2}=[diamond,draw=black,thick, inner sep=0pt ,minimum size=3.2mm]

\newtheorem{thm}{Theorem}
\newtheorem{lem}{Lemma}

\newtheorem{cor}{Corollary}

\newcommand{\qed}{$\Box$}
\newcommand{\smallqed}{{\tiny ($\Box$)}}

\newcommand{\barz}{\overline{z}}
\newcommand{\diam}{{\rm diam}}

\newcommand{\D}{{Dominator }}
\newcommand{\St}{{Staller }}

\newcommand{\gtd}{\gamma_{{\rm tg}}}
\newcommand{\gt}{\gamma_t}

\newcommand{\dstart}{\gamma_{{\rm tg}}}
\newcommand{\sstart}{\gamma_{{\rm tg}}^\prime}

\newenvironment{unnumbered}[1]{\trivlist \item [\hskip \labelsep {\bf
#1}]\ignorespaces\it}{\endtrivlist}

\newcommand{\cT}{{\cal T}}
\newcommand{\cF}{{\cal F}}

\newcommand{\proof}{\noindent\textbf{Proof. }}

\let\oldenumerate\enumerate
\renewcommand{\enumerate}{
  \oldenumerate
  \setlength{\itemsep}{0pt}
  \setlength{\parskip}{0pt}
  \setlength{\parsep}{0pt}
}

\begin{document}

\title{Trees with Equal Total Domination and \\ Game Total Domination Numbers}

\author{$^1$Michael A. Henning and $^2$Douglas F. Rall
\\ \\
$^1$Department of Pure and Applied Mathematics \\
University of Johannesburg \\
Auckland Park, 2006 South Africa\\
\small \tt Email: mahenning@uj.ac.za  \\
\\
$^2$Department of Mathematics \\
Furman University \\
Greenville, SC, USA\\
\small \tt Email: doug.rall@furman.edu}

\date{}
\maketitle

\begin{abstract}
In this paper, we continue the study of the total domination game in graphs introduced in [Graphs Combin. 31(5) (2015), 1453--1462], where the players Dominator and Staller alternately select vertices of $G$. Each vertex chosen must strictly increase the number of vertices totally dominated, where a vertex totally dominates another vertex if they are neighbors. This process eventually produces a total dominating set $S$ of $G$ in which every vertex is totally dominated by a vertex in $S$. Dominator wishes to minimize the number of vertices chosen, while Staller wishes to maximize it. The game total domination number, $\gamma_{\rm tg}(G)$, (respectively, Staller-start game total domination number, $\gamma_{\rm tg}'(G)$) of $G$ is the number of vertices chosen when Dominator (respectively, Staller) starts the game and both players play optimally. For general graphs $G$, sometimes $\gamma_{\rm tg}(G) > \gamma_{\rm tg}'(G)$. We show that if $G$ is a forest with no isolated vertex, then $\gamma_{\rm tg}(G) \le \gamma_{\rm tg}'(G)$. Using this result, we characterize the trees with equal total domination and game total domination number.
\end{abstract}


\newpage
\section{Introduction}

The domination game in graphs was first introduced by Bre\v{s}ar, Klav\v{z}ar, and Rall~\cite{BKR10} and extensively studied afterwards
in~\cite{BDKK14, brklra-2013, BKKR13, BKR13, Buj14, BuTu15, bujtas-2015+, doko-2013, HeLo-2015+, KWZ13, Kos, NaSiSo15} and elsewhere. Recently, the total version of the domination game was investigated in~\cite{hkr-2015}, where it was demonstrated that these two versions differ significantly. The total version has been studied in~\cite{BrHe16,BuHeTu16,BuHeTu17,dh-2016,hkr-2015+,HeRa16} and elsewhere. A vertex \emph{totally dominates} another vertex if they are neighbors. A \emph{total dominating set}, abbreviated TD-set, of a graph $G$ is a set $S$ of vertices such that every vertex of $G$ is totally dominated by a vertex in $S$. The \emph{total domination game} consists of two players called \emph{Dominator} and \emph{Staller}, who take turns choosing a vertex from~$G$. Each vertex chosen must totally dominate at least one vertex not totally dominated by the set of vertices previously chosen. Following the notation of~\cite{hkr-2015}, we call such a chosen vertex a \emph{legal move} or a \emph{playable vertex} in the total domination game. The game ends when the set of vertices chosen is a total dominating set in $G$. Thus we will  assume that all graphs under consideration in this paper have minimum degree at least~$1$.  Dominator's objective is to minimize the number of vertices chosen, while Staller's is to end the game with as many vertices chosen as possible.

The \emph{dominator-start total game} is the total domination game when Dominator starts the game, while the \emph{Staller-start total game} is the total domination game when \St starts the game.
The \emph{game total domination number}, $\gtd(G)$, of $G$ is the number of vertices chosen in the dominator-start total game when both players employ a strategy that achieves their objective. The number of vertices chosen in the Staller-start total game when both players employ a strategy that achieves their objective
is the \emph{Staller-start game total domination number}, $\gtd'(G)$, of $G$. Determining the exact value of $\gtd(G)$ and $\gtd'(G)$ is a challenging problem, and is currently known only for paths and cycles~\cite{dh-2016}.

Since the exact values of these invariants are very difficult to compute, we will often employ the so-called \emph{imagination strategy} when
it is required to show that the game total domination number of a tree and one of its subtrees differ by exactly (or by at most) some fixed amount.  This
method of proof was introduced in the initial paper~\cite{BKR10} on game domination.  Here it consists of both Dominator and Staller playing the total
domination game on one of the trees while Dominator ``imagines'' the game being played on the second tree.  Dominator chooses legal moves in the second
tree that are in keeping with his objective of minimizing the total number of vertices chosen there.  His moves in this second tree are then either
copied directly to the original tree (where the ``real'' game is being played by both players) or modified in some way so as to be legal moves
in the real game.  See~\cite{BKR10} for further explanation of this proof technique.

A \emph{partially total dominated graph} is a graph together with a declaration that some vertices are already totally dominated; that is, they need not be totally dominated in the rest of the game. Given a graph $G$ and a subset $S$ of vertices of $G$, we denote by $G|S$ the partially total dominated graph in which the vertices of $S$ in $G$ are already totally dominated. We use $\dstart(G|S)$ (resp. $\sstart(G|S)$) to denote the number of turns remaining in the game on $G|S$ under optimal play when Dominator (resp. Staller) has the next turn. In~\cite{hkr-2015}, the authors present a key lemma, named the \emph{Total Continuation Principle}.

\begin{lem}[``Total Continuation Principle'' -- \cite{hkr-2015}, Lemma~2.1]
Let $G$ be a graph and let $A,B \subseteq V(G)$.  If $B \subseteq A$, then $\dstart(G|A) \le \dstart(G|B)$ and $\sstart(G|A) \le \sstart(G|B)$.
\label{l:Cont}
\end{lem}

As a consequence of the Total Continuation Principle, when the total domination game is played on a partially total dominated graph $G$, the numbers $\gtd(G)$ and $\gtd'(G)$ can differ by at most~$1$.

\begin{cor}{\rm (\cite{hkr-2015})}
\label{c:diff1}
For every graph $G$ with no isolated vertex, we have $|\dstart(G) - \sstart(G)|\le 1$.
\end{cor}

\subsection{Notation}

For notation and graph theory terminology not defined herein, we in general follow~\cite{HeYe_book}. We denote
the \emph{degree} of a vertex $v$ in a graph $G$ by $d_G(v)$, or simply by $d(v)$ if the graph $G$ is clear
from the context. A \emph{degree-$k$ vertex} is a vertex of degree~$k$.
The minimum degree among the vertices of $G$ is denoted by $\delta(G)$.
A vertex of degree~$1$ is called a \emph{leaf} and its neighbor a \emph{support vertex}. A \emph{strong support vertex} is a support vertex with at least two leaf neighbors. A star is a tree with at most one vertex of degree~$2$ or more. A \emph{subdivided star} is the tree obtained from a star on at least three vertices by subdividing every edge exactly once.
The \emph{open neighborhood} of a vertex $v \in V(G)$ is $N_G(v) = \{u \in V(G) \, | \, uv \in E(G) \}$ and the degree of $v$ is
$d_G(v) = |N_G(v)|$. The \emph{closed neighborhood of $v$} is $N_G[v] = \{v\} \cup N_G(v)$.

For a set $S \subseteq V(G)$, we let $G[S]$ denote the subgraph induced by $S$. The graph obtained from $G$ by deleting the vertices
in $S$ and all edges incident with vertices in $S$ is denoted by $G-S$. If $S = \{v\}$, we also denote $G-S$ simply by $G-v$.

If $X$ and $Y$ are subsets of vertices in a graph $G$, then the set $X$ \emph{totally dominates} the set
$Y$ in $G$ if every vertex of $Y$ is adjacent to at least one vertex of $X$. In particular, if $X$ totally
dominates the vertex set of $G$, then $X$ is a TD-set in $G$. The cardinality of a smallest TD-set in $G$ is the \emph{total domination number} of $G$ and is denoted $\gt(G)$.
A TD-set of $G$ of cardinality $\gt(G)$ is called a
$\gt(G)$-\emph{set}.  Since an isolated vertex in a graph cannot be totally dominated by definition, all graphs considered will be without isolated vertices.
For more information on total domination in graphs see the recent book~\cite{HeYe_book}.
We use the standard notation $[k] = \{1,\ldots,k\}$.

A \emph{rooted tree} $T$ distinguishes one vertex $r$ called the \emph{root}. For each vertex $v \ne r$ of $T$, the \emph{parent} of $v$ is the
neighbor of $v$ on the unique $(r,v)$-path, while a \emph{child} of $v$ is any other neighbor of $v$. We denote all the children of a vertex $v$
by $C(v)$. A \emph{descendant} of $v$ is a vertex $u \ne v$ such that the unique $(r,u)$-path contains $v$. Thus, every child of $v$ is a descendant of~$v$. An \emph{ancestor} of $v$ is a vertex $u \ne v$ that belongs to the $(r,v)$-path in $T$. In particular, the parent of $v$ is an ancestor of $v$. The \emph{grandparent} of $v$ is the ancestor of $v$ at distance~$2$ from $v$. A \emph{grandchild} of $v$ is the descendant of $v$ at distance~$2$ from $v$. A path on $n$ vertices is denoted by $P_n$.

Let $G$ be a partially total dominated graph and let $v$ be a vertex of $G$. If $v$ is not totally dominated in $G$, we call the vertex
$v$ \emph{totally undominated} in $G$. We let $G_v$ denote the partially total dominated graph obtained from $G$ by totally dominating $N(v)$.
If the vertex $v$ is totally dominated in $G$, then we let $G^v$ denote the partially total dominated graph obtained from $G$ by removing $v$
from the set of totally dominated vertices. We note that $G^v$ and $G$ are identical except that $v$ is totally dominated in $G$ but not in $G^v$.

\section{Main Result}

As remarked by Bre\v{s}ar, Klav\v zar, Ko\v{s}mrlj, and Rall~\cite{BKKR13}, ``the domination game is very non-trivial even when played on trees." In this paper we prove the following result.

\begin{thm}
If $F$ is a partially total dominated forest with no isolated vertex, then $\dstart(F) \le \sstart(F)$.
\label{t:main1}
\end{thm}

We remark that Theorem~\ref{t:main1} is not true for general graphs. For example, $\dstart(C_8) = 5 = \sstart(C_8) + 1$.
As an application of Theorem~\ref{t:main1}, we prove our main result, which is a characterization of trees with equal total domination and game total domination number. For this purpose, we construct a family~$\cF$ of trees with equal total domination and game total domination number in Section~\ref{S:cF}. A \emph{nontrivial} tree is a tree of order at least~$2$.

\begin{thm}
Let $T$ be a nontrivial tree. Then, $\gt(T) = \dstart(T)$ if and only if $T \in \cF$.
\label{t:main2}
\end{thm}

We proceed as follows. We present a proof of our first main result, namely Theorem~\ref{t:main1}, in Section~\ref{S:proof1}. Thereafter, we present a series of preliminary lemmas in Section~\ref{S:prelim}, before giving a proof of our second main result, namely Theorem~\ref{t:main2}, in Section~\ref{S:proof2}.

\section{Proof of Theorem~\ref{t:main1}}
\label{S:proof1}

In this section we present a proof of Theorem~\ref{t:main1}. For this purpose, we introduce some additional notation. We shall need the following property of partially total dominated forests.

\begin{lem}
\label{lem:main1}
Given $\ell \ge 2$, assume $\dstart(F) \le \sstart(F)$ for all partially total dominated forests $F$ such that $\dstart(F) \le \ell$. If $F$ is such a forest that contains two neighbors $x$ and $y$ such that the vertex~$x$ and every vertex at distance~$1$ and~$2$ from $x$, except possibly for $y$, is totally dominated in $F$, then $\dstart(F^{x}) > \dstart(F)$ and $\sstart(F^{x}) > \sstart(F)$.
\end{lem}
\proof Fixing $\ell$, we use induction on the number of vertices in $F$ that are totally undominated to prove that $\dstart(F^{x}) > \dstart(F)$ and $\sstart(F^{x}) > \sstart(F)$. If there are no totally undominated vertices in $F$, then $\dstart(F) = 0$ since every vertex is totally dominated in $F$. Further, $\dstart(F^{x}) = \sstart(F^{x}) = 1$, since $x$ is the only vertex that is totally undominated in $F^{x}$. Thus, trivially $\dstart(F^{x}) > \dstart(F)$ and $\sstart(F^{x}) > \sstart(F)$. This establishes the base case. Let $k \ge 1$ and suppose that the desired result holds if there are fewer than~$k$ totally undominated vertices in $F$. Let $F$ have $k$ totally undominated vertices.

To prove $\dstart(F^{x}) > \dstart(F)$, let $v$ be an optimal first move in the Dominator-start game played in $F^{x}$. We may assume that $\dstart(F^{x}) \le \ell$, for otherwise $\dstart(F^{x}) \ge \ell+1 > \ell \ge \dstart(F)$.
If $v$ is a neighbor of $x$, then $x$ is the only new vertex totally dominated by~$v$ in $F^{x}$ since, by supposition, every vertex at distance~$1$ and~$2$ from $x$, except possibly for $y$, is totally dominated in $F$ and therefore also in $F^{x}$. Thus, $\dstart(F^{x}) = 1 + \sstart(F)$. By assumption, $\sstart(F) \ge \dstart(F)$, implying that  $\dstart(F^{x}) > \dstart(F)$. Hence, we may assume that $v$ is not a neighbor of $x$, for otherwise the desired result holds.

By the choice of $v$, we have $\dstart(F^{x})=1+\sstart(F^{x}_v)$, and by the Total Continuation Principle, $\dstart(F_v) \le \dstart(F) \le \ell$.
By assumption the forest $F$ has exactly $k$ totally undominated vertices, and since $v$ is a legal move in $F^{x}$ and $v$ is not a neighbor of $x$, it follows that
$F_v$ has fewer than $k$ totally undominated vertices.  Furthermore, in $F_v$ the vertex $x$  and every vertex at distance
1 or 2 from $x$, except possible for $y$, is totally dominated.  Therefore, we can apply the induction hypothesis to $F_v$.  We get
$\sstart(F^x_v)> \sstart(F_v)$. We note that the vertex~$v$ may not be an optimal first move for \D in $F$. Thus, since \D does at least as well in $F$ by playing optimally as by playing $v$ first, $\dstart(F) \le 1 + \sstart(F_v)$. These observations imply that
\[
\dstart(F^{x}) = 1 + \sstart(F^{x}_v) >  1 + \sstart(F_v) \ge 1 + (\dstart(F) - 1) = \dstart(F).
\]

To prove $\sstart(F^{x}) > \sstart(F)$, let $v$ be an optimal first move in the Staller-start total game played in $F$. By optimality of $v$, we
have $\sstart(F) = 1 + \dstart(F_v)$. Since $F$ is a forest, and since the vertex $x$ and every vertex at distance~$1$ and~$2$ from~$x$, except
possibly for $y$, is totally dominated in $F$, we note that $v$ is not a neighbor of~$x$. Since $v$ is not adjacent to $x$, we have $(F_v)^x=(F^x)_v$ and
we denote this forest simply by $F_v^x$.  Since $v$ is a legal move in $F$, there are fewer than $k$ totally undominated vertices in the partially total dominated forest $F_v$. Also, it follows from the Total Continuation Principle that $\dstart(F_v) \le \dstart(F) \le \ell$. Applying the
inductive hypothesis to $F_v$, we have $\dstart(F^{x}_v) > \dstart(F_v)$. We note that the vertex~$v$ may not be an optimal first move for \St in~$F^{x}$. Thus, since \St does at least as well in $F^{x}$ by playing optimally as by playing $v$ first, $\sstart(F^{x}) \ge 1 + \dstart(F^{x}_v)$. These observations
imply that
\[
\sstart(F^{x}) \ge 1 + \dstart(F^{x}_v) > 1 + \dstart(F_v) = 1 + (\sstart(F) - 1) = \sstart(F). \hspace*{0.5cm} \Box
\]

\medskip
We are now in a position to prove Theorem~\ref{t:main1}. Recall its statement. We remark that our proof of Theorem~\ref{t:main1} employs some of the key ideas from a proof of an analogous result for the ordinary game domination number due to Kinnersley, West, and Zamani~\cite{KWZ13}.

\noindent \textbf{Theorem~\ref{t:main1}}. \emph{If $F$ is a partially total dominated forest with no isolated vertex, then $\dstart(F) \le \sstart(F)$.
}

\medskip
\proof  By Corollary~\ref{c:diff1}, every forest $F$ with no isolated vertex satisfies $|\dstart(F) - \sstart(F)| \le 1$. Hence, it suffices for us to
prove that for all $k \ge 1$, $\dstart(F) = k$ and $\sstart(F) = k - 1$ cannot both hold. If $\dstart(F) = 1$, then the forest $F$ contains some totally
undominated vertices, and so $\sstart(F) \ge 1$. If $\sstart(F) = 1$, then every legal move completes the game, and so $\dstart(F) = 1$. Thus, if $k \in \{1,2\}$,
then $\dstart(F) = k$ and $\sstart(F) = k - 1$ cannot both hold. This establishes the base case. Let $k \ge 3$ and assume that if $F'$ is a partially total dominated forest with no isolated vertex satisfying $\dstart(F') < k$, then $\dstart(F') \le \sstart(F')$. Let $F$ be a partially total dominated forest with
no isolated vertex satisfying $\dstart(F) = k$, and suppose, to the contrary, that $\sstart(F) = k - 1$. We proceed further with the following claim.

\begin{unnumbered}{Claim~\ref{t:main1}.1}
For every playable vertex $v$ in $F$, we have $\dstart(F_v) = k-2$ and $\sstart(F_v) = k-1$.
\end{unnumbered}
\ \\[-2em]
\noindent\textbf{Proof.}
Let $v$ be any playable vertex in $F$, and so at least one neighbor of $v$ is totally undominated in $F$. Since $\sstart(F) = k - 1$,
after any first move of Staller, \D can complete the game by forcing at most~$k-2$ additional moves to be played. In particular, if \St plays the vertex~$v$
as her first move, then $\dstart(F_v) \le k-2$. Since $\dstart(F) = k$, after any first move of Dominator, \St can complete the game by forcing at least~$k-1$
additional moves to be played. In particular, if \D plays the vertex~$v$ as his first move, then $\sstart(F_v) \ge k-1$. Thus, for every playable vertex $v$
in $F$, we have $\dstart(F_v) \le k-2$ and $\sstart(F_v) \ge k-1$. By Corollary~\ref{c:diff1}, it now follows that $\dstart(F_v) = k-2$ and
$\sstart(F_v) = k-1$.~\smallqed

\medskip
We now return to the proof of Theorem~\ref{t:main1}. If every component of $F$ contains at most one totally undominated vertex, then
$\dstart(F) = \sstart(F)$, a contradiction. Therefore, there is a component $C$ of $F$ containing at least two totally undominated vertices.
We now root the component $C$ at an arbitrary vertex $r$ of $C$, and let $x$ be a totally undominated vertex at maximum distance
from~$r$ in $C$. By our choice of the component $C$, we note that $x \ne r$. Let $y$ be the parent of $x$ in the rooted tree $C$.
We now consider the partially total dominated forest $F_y$. By our choice of the
vertex~$x$, all descendants of $x$ are totally dominated in $F$. Further, $x$ and $y$ are neighbors in $F_y$ such that the vertex~$x$ and every vertex
at distance~$1$ and~$2$ from $x$, except possibly for $y$, is totally dominated in $F_y$. This implies that no neighbor of $x$ is a legal move in $F_y$.
Let $v$ be an optimal first move in the Staller-start total game on $F_y$. Thus, $v \notin N(x)$.

We now consider the partially total dominated forest $(F_y)_v$  obtained from $F_y$ by totally dominating $N(v)$. We denote this graph simply
by $F_{y,v}$. By Claim~\ref{t:main1}.1 and by the optimality of the vertex $v$, we have $k - 1 = \sstart(F_y) = 1 + \dstart(F_{y,v})$. Thus,
$\dstart(F_{y,v}) = k-2$.
We now consider the partially total dominated forest $F_{y,v}^x$ obtained from $F_{y,v}$ by removing $x$ from the set of totally dominated
vertices in $F_{y,v}$.  We note that $F_{y,v}^x$ and $F_{y,v}$ are identical except that $x$ is totally dominated in $F_{y,v}$ but not in $F_{y,v}^x$.
Further, we note that every vertex at distance~$1$ and~$2$ from $x$, except possibly for $y$, is totally dominated in both $F_{y,v}$ and $F_{y,v}^x$.
Since $\dstart(F_{y,v})=k-2$, it follows from the inductive hypothesis that $\dstart(F_{y,v}) \le \sstart(F_{y,v})$.
Applying Lemma~\ref{lem:main1} to the partially total dominated forest $F_{y,v}$, we have
\[
\dstart(F_{y,v}^x) > \dstart(F_{y,v}) \hspace*{0.5cm} \mbox{and} \hspace*{0.5cm}  \sstart(F_{y,v}^x) > \sstart(F_{y,v}).
\]

By Claim~\ref{t:main1}.1, we note that $\dstart(F_{v}) = k-2$. By the Total Continuation Principle, we have $\dstart(F_{v}) \ge \dstart(F_{y,v}^x)$. As
observed earlier, $\dstart(F_{y,v}) = k-2$. Thus,
\[
k-2 = \dstart(F_{v}) \ge \dstart(F_{y,v}^x) > \dstart(F_{y,v}) = k-2,
\]
a contradiction. Therefore our supposition that $\sstart(F) = k - 1$ is false, implying that if $\dstart(F) = k$, then $\sstart(F) \ge k$. This
completes the proof of Theorem~\ref{t:main1}.~\qed

\section{Preliminary Lemmas}
\label{S:prelim}

In this section, we present some preliminary lemmas.
By the Total Continuation Principle, it is never in Dominator's best interests to play a leaf that belongs to a component in the partially total dominated forest that is not a star, since in this case \D can always do at least as well by playing a non-leaf neighbor of a support vertex instead of one of its leaf-neighbors.

We begin with the following properties of trees with equal total domination and game total domination numbers. The first two results hold for graphs in general.

\begin{lem}
\label{lem1}
Let $G$ be a graph with no isolated vertex satisfying $\gt(G) = \dstart(G)$. Every (legal) move that Staller can play on each of her turns in the
dominator-start total game played in $G$ is an optimal move for her in the sense that if \D plays optimally, he always finishes the game in exactly
$\dstart(G)$ moves and no fewer, whatever choice of moves \St makes.
\end{lem}
\proof If \D plays optimally in the dominator-start total game and \St plays any legal move on each of her turns, then \D can guarantee that the
game requires at most~$\dstart(G)$ moves, and possibly fewer if \St does not play optimally. However, since the graph $G$ satisfies $\gt(G) = \dstart(G)$,
the resulting set of played vertices in the dominator-start total game is a minimum TD-set in $G$. This implies that Staller's moves have
no bearing on the outcome of the game, in that \D cannot finish the game in fewer than~$\dstart(G)$ moves, whatever choice of (legal) moves \St plays.
Thus, in this case, every move that \St makes is an optimal move.~\qed

\medskip
During the course of the total domination game, we say that \D can \emph{block} a move $v$ of \St if he can play a vertex that results in all neighbors of $v$ totally dominated, implying that $v$ is not a legal move in the remaining part of the game.

\begin{lem}
\label{lem2}
If $G$ is a graph with no isolated vertex satisfying $\gt(G) = \dstart(G)$, then no two degree-$1$ vertices are at distance~$3$ apart in $G$.
\end{lem}
\proof Suppose, to the contrary, that $u$ and $v$ are two degree-$1$ vertices at distance~$3$ apart in $G$. No first (optimal) move of \D can block both $u$ and $v$. Hence on Staller's first move, she can play $u$ or $v$. However, neither $u$ nor $v$ belong to a minimum TD-set in $G$, implying that $\gt(G) < \dstart(G)$, a contradiction.~\qed

\begin{figure}
\begin{center}
\begin{tikzpicture}[scale=.8,style=thick,x=1cm,y=1cm]
\def\vr{2.5pt} 
\path (1.25,0.75) coordinate (x1);
\path (1.75,0.75) coordinate (x2);
\path (1.5,0.5) coordinate (v);
\path (1,0) coordinate (w);
\path (2,0) coordinate (u);
%
\draw (v) -- (x1);
\draw (v) -- (x2);
\draw (v) -- (w);
\draw (v) -- (u);
\draw (v) [fill=black] circle (\vr);
\draw (w) [fill=black] circle (\vr);
\draw (u) [fill=black] circle (\vr);
\draw[anchor = south] (v) node {$v$};
\draw[anchor = east] (w) node {$u_1$};
\draw[anchor = west] (u) node {$u_2$};
\draw (1.5,-1) node {(a) A $(1,1)$-structure};

\path (9.15,1.85) coordinate (x1);
\path (9.85,1.85) coordinate (x2);
\path (9.5,1.5) coordinate (v);
\path (8.8,0.75) coordinate (v1);
\path (8.8,0) coordinate (v2);
\path (10.2,0.75) coordinate (w1);
\path (10.2,0) coordinate (w2);
%
\draw (v) -- (x1);
\draw (v) -- (x2);
\draw (v) -- (v1);
\draw (v) -- (w1);
\draw (v1) -- (v2);
\draw (w1) -- (w2);
\draw (v) [fill=black] circle (\vr);
\draw (v1) [fill=black] circle (\vr);
\draw (v2) [fill=black] circle (\vr);
\draw (w1) [fill=black] circle (\vr);
\draw (w2) [fill=black] circle (\vr);
\draw[anchor = south] (v) node {$v$};
\draw[anchor = east] (v1) node {$v_1$};
\draw[anchor = east] (v2) node {$u_1$};
\draw[anchor = west] (w1) node {$v_2$};
\draw[anchor = west] (w2) node {$u_2$};
\draw (9.5,-1) node {(b) A $(2,2)$-structure};
\end{tikzpicture}
\end{center}
\vskip -0.6 cm
\caption{A $(1,1)$-structure and a $(2,2)$-structure with pivot~$v$.} \label{f:S11}
\end{figure}
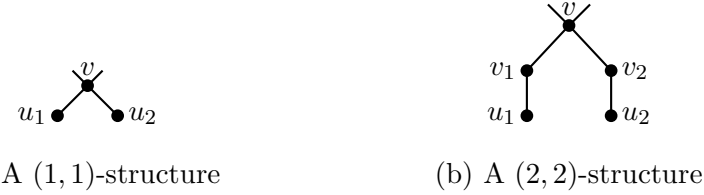

\medskip
We introduce next some additional notation. For $r \ge 1$, we define an $(\ell_1,\ldots, \ell_r)$-structure pivoted at a vertex~$v$ in a tree $T$ to be
$r$ paths $Q_1, \ldots, Q_r$ of lengths~$\ell_1, \ldots,
\ell_r$, respectively, emanating from $v$ in $T$, such that the degree of every vertex on these $r$ paths, except possibly for the vertex~$v$, is the
same as its degree in $T$. We call the vertex $v$ the \emph{pivot vertex} of the associated $(\ell_1,\ldots, \ell_r)$-structure, which we denote by $T_v$.
A $(1,1)$-structure and a $(2,2)$-structure are illustrated in Figure~\ref{f:S11}(a) and~\ref{f:S11}(b), respectively.

\begin{lem}
\label{lem3}
Let $T$ be a tree satisfying $\gt(T) = \dstart(T)$ that contains a $(1,1)$-structure, $T_v$, with pivot vertex~$v$. If $T'$ is the tree obtained from $T$
by deleting a leaf vertex in $T_v$, then $\gt(T') = \dstart(T')$.
\end{lem}
\proof Let $u_1$ and $u_2$ be the two leaf neighbors of $v$ in $T_v$. Renaming vertices if necessary, we may assume that $T' = T - u_1$. The vertex $v$ is
a support vertex in $T'$, and therefore belongs to every TD-set of $T'$. In particular, every $\gt(T')$-set  is a TD-set of $T$,
and so $\gt(T) \le \gt(T')$.
Conversely, there is a $\gt(T)$-set that does not contain $u_1$.  This set is  also
a TD-set of $T'$, and so $\gt(T') \le \gt(T)$. Consequently, $\gt(T') = \gt(T)$.
We show next that $\dstart(T') \le \dstart(T)$. Consider the dominator-start total game played in $T'$. \D plays an imaginary strategy in $T$, although the
real game is played in $T'$. Each of Staller's moves is played in $T'$. However, \D imagines the game to be played in $T$ and on each of his turns, he
considers an optimal move that would be played in~$T$. Suppose that the leaf $u_1$ is played in the imaginary strategy in $T$. The only vertex totally
dominated by $u_1$ is its neighbor~$v$. Immediately before $u_1$ is played, it is an optimal move for Dominator. Since both leaves $u_1$ and $u_2$ share a
common neighbor, namely~$v$, in $T$, the leaf $u_2$ is therefore also an optimal move for Dominator. Thus, renaming vertices, if necessary, we may assume
that if $u_1$ or $u_2$ is played, then the leaf $u_2$ is played instead of the leaf~$u_1$ in the imaginary game in~$T$. With this assumption, every (optimal) move played by \D in the imaginary game in $T$ is a legal move in $T'$. Clearly, every move of \St in $T'$ is a playable vertex in~$T$.

\D now imagines each of Staller's moves to be played in $T$, and considers an optimal move in $T$ that he would play in response to her move. By our assumption
that every (optimal) move played by \D in the imaginary game in $T$ is a legal move in $T'$, \D plays an optimal move in $T$ on each of his moves in $T'$.
By Lemma~\ref{lem1}, each of Staller's moves is an optimal move in the imaginary game played in $T$. Thus, Dominator's strategy of playing the imaginary game
guarantees that the (real) game in $T'$ is finished in exactly~$\dstart(T)$ moves. Thus, $\dstart(T') \le \dstart(T)$. As observed earlier, $\gt(T') = \gt(T)$.
By assumption, $\gt(T) = \dstart(T)$. Therefore,
$\dstart(T) = \gt(T) = \gt(T') \le \dstart(T') \le \dstart(T)$. Consequently, we must have equality throughout the above inequality chain. In particular,
$\gt(T') = \dstart(T')$. This completes the proof of Lemma~\ref{lem3}.~\qed

\medskip
An analogous proof to that of Lemma~\ref{lem3} establishes the following result.

\begin{lem}
\label{lem_strong}
If $T$ is obtained from a nontrivial tree $T'$ by adding a pendant edge to a support vertex of $T'$, then $\dstart(T) = \dstart(T')$.
\end{lem}

As an immediate consequence of Lemma~\ref{lem2}, we have the following result.

\begin{lem}
\label{lem4}
If $T$ is a tree satisfying $\gt(T) = \dstart(T)$, then $T$ contains no $(2,1)$-structure.
\end{lem}

In the subsequent lemmas, we assume throughout that $T$ is a tree satisfying $\gt(T) = \dstart(T)$ and $T$ contains no $(1,1)$-structure (that is, $T$ has no strong support vertex). By Lemma~\ref{lem4}, the tree $T$ contains no $(2,1)$-structure. Further, we assume that $T$ is rooted at an optimal first move, $r$ say, of Dominator. Thus, \D plays the vertex $d_1 = r$ as his first move in the dominator-start total game played in $T$. Let $u$ be a vertex at maximum distance from the root $r$ in $T$. Necessarily, $u$ is a leaf. Let $v$ be the parent of $u$ and let $w$ the parent of $v$. Further, if $r \ne w$, let $x$ the parent of $w$, and if $r \ne x$, let $y$ the parent of $x$.

\begin{lem}
If the leaf $u$ belongs to a $(2,2)$-structure, $T_w$, with pivot vertex~$w$, then $T$ is a subdivided star.
\label{lem5}
\end{lem}
\proof Let $T_w$ be the path $uvwv'u'$, where $v$ and $v'$ have degree~$2$ in $T$ and $u$ and $u'$ are leaves in $T$. We note that every $\gamma_t(T)$-set contains the three vertices $v$, $v'$ and $w$, and therefore $u$ and $u'$ belong to no $\gamma_t(T)$-set. If the root $r$ is not the vertex $w$, then \St can play as her first move the vertex $u$, implying that $\gt(T) < \dstart(T)$, a contradiction. Hence, $r = w$. Since $T$ has no $(2,1)$-structure by Lemma~\ref{lem4}, this implies that $T$ is a subdivided star.~\qed

\medskip
In what follows, we may assume that a vertex at maximum distance from the root $r$ in $T$ does not belong to a $(2,2)$-structure.

\begin{figure}[htb]
\begin{center}
\begin{tikzpicture}[scale=.8,style=thick,x=1cm,y=1cm]
\def\vr{2.5pt} 
\path (0,0) coordinate (u);
\path (0,0.75) coordinate (v);
\path (0,1.5) coordinate (w);
\path (0.7,2.25) coordinate (x);
\path (0.7,3.25) coordinate (r);
\path (1.4,1.5) coordinate (w1);
\path (0.7,2.7) coordinate (x1);
\path (1.5,1.75) coordinate (x2);
\path (1.6,1.95) coordinate (x3);
\path (1.5,2.75) coordinate (r2);
\path (-0.1,2.75) coordinate (r3);
\draw (x) -- (r);
\draw (x) -- (x2);
\draw (x) -- (x3);
\draw (x) -- (w1);
\draw (x) -- (w);
\draw (w) -- (v);
\draw (v) -- (u);
\draw (r) -- (r2);
\draw (r) -- (r3);
\draw (u) [fill=black] circle (\vr);
\draw (v) [fill=black] circle (\vr);
\draw (w) [fill=black] circle (\vr);
\draw (w1) [fill=black] circle (\vr);
\draw (x) [fill=black] circle (\vr);
\draw (r) [fill=black] circle (\vr);
\draw[anchor = south] (r) node {$r$};
\draw[anchor = east] (x) node {$x$};
\draw[anchor = east] (u) node {$u$};
\draw[anchor = east] (v) node {$v$};
\draw[anchor = east] (w) node {$w$};
\draw[anchor = west] (w1) node {$w'$};
\end{tikzpicture}
\end{center}
\vskip -0.6 cm
\caption{The $(3,1)$-structure $T_x$ with pivot~$x$.} \label{f:S31}
\end{figure}
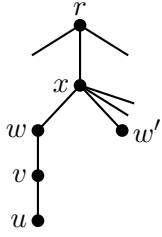

\begin{lem}
If the leaf $u$ belongs to a $(3,1)$-structure, $T_x$, with pivot vertex~$x$ and $T' = T - \{u,v,w\}$, then $r = y$ and $\gt(T') = \dstart(T')$.
\label{lem6}
\end{lem}
\proof Let $w'$ be the leaf-neighbor of $x$ in the $(3,1)$-structure $T_x$. Thus, $T_x$ is the path $uvwxw'$, where $v$ and $w$ have degree~$2$ in $T$ and $u$ and $w'$ are leaves. We note that the support vertices $v$ and $x$ belong to every $\gamma_t(T)$-set. If $w'$ belongs to some $\gamma_t(T)$-set $S$, then $w \notin S$, and so $u \in S$. But then $(S \setminus \{u,w'\}) \cup \{w\}$ is a TD-set of $T$ of size less than~$|S|$, a contradiction. Thus, $w'$ belongs to no $\gamma_t(T)$-set. If the root $r$ is not the parent of $x$, then \St can play as her first move the vertex $w'$, implying that $\gt(T) < \dstart(T)$, a contradiction. Hence, $r$ is the parent of $x$; that is, $r = y$. The $(3,1)$-structure $T_x$ is illustrated in Figure~\ref{f:S31}.

We now consider the tree $T' = T - \{u,v,w\}$.  As in the proof of Lemma~\ref{lem5}, \D plays an imaginary game in $T$, although the real
game is played in $T'$. By Lemma~\ref{lem1}, each of Staller's moves in the game played in $T'$ is an optimal move in the imaginary game played in $T$. Each of Staller's moves in the real game is played in $T'$. However, \D imagines the game to be played in $T$ and on each of his turns, he considers an optimal move that would be played in~$T$. Since \D plays the vertex $d_1 = r$ as his first move, and since $d_1$ totally dominates the vertex $x$, we note that every move of \St in the real game played in $T'$ is a legal move in the imaginary game played in $T$.

We show first that $\dstart(T') \le \dstart(T) - 2$. If every (optimal) move of \D in response to each of Staller's moves is a vertex in $T'$, then the game finishes in at most~$\dstart(T) - 2$ moves, since in the imaginary game played in $T$ at least two further moves are needed (in order to totally dominate $u$ and $v$). Thus, in this case, \D has a strategy to complete the game in $T'$ in at most~$\dstart(T) - 2$ moves. Hence, we may assume that during the imaginary game played in $T$, \D plays a vertex not in the real game $T'$. As observed earlier, by the Total Continuation Principle, it is never in Dominator's best interests to play a leaf. Thus, the first vertex played by \D in the imaginary game that is not in $T'$ is the vertex $z$, where $z \in \{v,w\}$. Suppose that this occurs after Staller's $k$th move; that is, after $2k$ vertices have been played. Let $\barz = \{v,w\} \setminus \{z\}$. \D imagines that \St immediately responds by playing $\barz$ as her $(k+1)$st move, which by Lemma~\ref{lem1}, is an optimal move in the imaginary game played in $T$. We note that neither $v$ nor $w$ totally dominate any new vertex in the real game played in $T'$. As his $(k+1)$st move in the real game, \D then plays an optimal move which he would play in the imaginary game in response to Staller's imagined move~$\barz$. From here onwards, all remaining moves of \D (in the imaginary game played in $T$) are played in the real game $T'$. Thus, once again, \D has a strategy to complete the game in $T'$ in at most~$\dstart(T) - 2$ moves. Therefore, $\dstart(T') \le \dstart(T) - 2$.

Every $\gt(T')$-set can be extended to a TD-set of $T$ by adding to it the vertices $v$ and $w$, implying that $\gt(T) \le \gt(T') + 2$. By assumption, $\gt(T) = \dstart(T)$. As observed earlier, $\dstart(T') \le \dstart(T) - 2$. Therefore, since
$\gt(G) \le \dstart(G)$ holds for every (isolate-free) graph $G$, we have
\[
\dstart(T) = \gt(T) \le \gt(T') + 2 \le \dstart(T') + 2 \le \dstart(T).
\]
Consequently, we must have equality throughout the above inequality chain. In particular, $\gt(T') = \dstart(T')$. This completes the proof of
Lemma~\ref{lem6}.~\qed

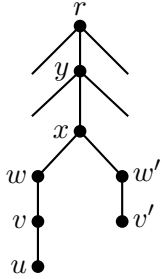
\begin{figure}[htb]
\begin{center}
\begin{tikzpicture}[scale=.8,style=thick,x=1cm,y=1cm]
\def\vr{2.5pt} 
\path (0,0) coordinate (u);
\path (0,0.75) coordinate (v);
\path (0,1.5) coordinate (w);
\path (0.7,2.25) coordinate (x);
\path (0.7,3.25) coordinate (y);
\path (0.7,4) coordinate (r);
\path (1.4,1.5) coordinate (w1);
\path (1.4,1.6) coordinate (wp);
\path (1.4,0.75) coordinate (v1);
\path (1.4,0.85) coordinate (vp);
\path (1.5,2.75) coordinate (y2);
\path (-0.1,2.75) coordinate (y3);
\path (1.5,3.2) coordinate (r2);
\path (-0.1,3.2) coordinate (r3);
\path (1.5,2.5) coordinate (y2);
\path (-0.1,2.5) coordinate (y3);
\draw (x) -- (y);
\draw (y) -- (r);
\draw (y) -- (y2);
\draw (y) -- (y3);
\draw (r) -- (r2);
\draw (r) -- (r3);
\draw (x) -- (w1);
\draw (w1) -- (v1);
\draw (x) -- (w);
\draw (w) -- (v);
\draw (v) -- (u);
\draw (u) [fill=black] circle (\vr);
\draw (y) [fill=black] circle (\vr);
\draw (v) [fill=black] circle (\vr);
\draw (w) [fill=black] circle (\vr);
\draw (w1) [fill=black] circle (\vr);
\draw (v1) [fill=black] circle (\vr);
\draw (x) [fill=black] circle (\vr);
\draw (r) [fill=black] circle (\vr);
\draw[anchor = south] (r) node {$r$};
\draw[anchor = east] (x) node {$x$};
\draw[anchor = east] (y) node {$y$};
\draw[anchor = east] (u) node {$u$};
\draw[anchor = east] (v) node {$v$};
\draw[anchor = east] (w) node {$w$};
\draw[anchor = west] (wp) node {$w'$};
\draw[anchor = west] (vp) node {$v'$};
\end{tikzpicture}
\end{center}
\vskip -0.6 cm
\caption{The $(3,2)$-structure $T_x$ with pivot~$x$.} \label{f:S32}
\end{figure}

\begin{lem}
\label{lem7}
If the leaf $u$ belongs to a $(3,2)$-structure, $T_x$, with pivot vertex~$x$, then one of the following holds. \\[-27pt]
\begin{enumerate}
\item The root $r = x$, and $T$ is obtained from a star with $k_1 + k_2$ leaves, where $k_1,k_2 \ge 1$, by subdividing $k_1$ edges once and $k_2$ edges twice.
\item The root $r \ne x$ and the vertex $x$ has degree~$3$. Further, the root $r$ is the parent of $y$, and if $T' = T - V(T_x)$, then $\gt(T') = \dstart(T')$.
\end{enumerate}
\end{lem}
\proof Let $T_x$ be the path $uvwxw'v'$. Thus, $v$, $w$ and $w'$ have degree~$2$ in $T$ and $u$ and $v'$ are leaves. If the root $r = x$, then Part~(a) follows from Lemma~\ref{lem4} and our assumption that there is no strong support vertex. Hence, we may assume that the root $r \ne x$. Thus, \St can play as her first move the leaf $v'$; that is, $s_1 = v'$.  If \D or \St can play the vertex $x$ on any move during the remainder of the game,  then this would imply that both $x$ and $v'$ belong to a $\gamma_t(T)$-set, a contradiction. Hence, \D is forced to play the vertex $d_2 = v$ as his second move in order to block the vertex $x$ from being played. Further, after \D plays the vertex $v$, the vertex $x$ must not be playable, implying that $x$ has degree~$3$ in $T$ and that the root $r$ is the grandparent of $x$ (or, equivalently, the parent of $y$). The tree $T$, with $(3,2)$-structure $T_x$, is illustrated in Figure~\ref{f:S32}.

We now consider the tree $T' = T - V(T_x)$ and show that $\gt(T') = \dstart(T')$. As before, \D plays an imaginary game in $T$, although the real game is played in $T'$. Since \D plays the vertex $d_1 = r$ as his first move, and since $d_1$ totally dominates the vertex $y$, we note that every move of \St in the real game played in $T'$ is a legal move in the imaginary game played in $T$.

We show first that $\dstart(T') \le \dstart(T) - 4$. If every (optimal) move of \D in response to each of Staller's moves is a vertex in $T'$, then the game finishes in at most~$\dstart(T) - 4$ moves, since in the imaginary game played in $T$ at least four further moves are needed (in order to totally dominate the four vertices $u$, $v$, $v'$ and $w'$). Thus, in this case, \D has a strategy to complete the game in $T'$ in at most~$\dstart(T) - 4$ moves. Hence, we may assume that during the imaginary game played in $T$, \D plays a vertex not in the real game $T'$. By the Total Continuation Principle, the first vertex played by \D in the imaginary game that is not in $T'$ is the vertex $z$, where $z \in \{v,w,w',x\}$. Suppose that this occurs after Staller's $k$th move; that is, after $2k$ vertices have been played. If $z = v$ or $z = w$, then \D imagines that \St immediately responds by playing $w$ or $v$, respectively, as her $(k+1)$st move. If $z = w'$ or $z = x$, then \D imagines that \St immediately responds by playing $v'$ or $w'$, respectively, as her $(k+1)$st move. In all four cases, Staller's move is an optimal move in the imaginary game played in $T$ by Lemma~\ref{lem1}.

If now Dominator's $(k+2)$nd move in the imaginary game belongs to $V(T')$, then he plays this move as his $(k+1)$st move in the real game. If every subsequent move of \D in response to a move of Staller is a vertex in $T'$, then the game finishes in at most~$\dstart(T) - 4$ moves, since in the imaginary game played in $T$ at least two further moves are needed (in order to totally dominate all vertices in $V(T_x)$). Hence, we may assume that as his $\ell$th move in the real game, where $\ell \ge k+2$, \D plays a vertex not in the real game $T'$. Analogously as before, \D imagines that \St immediately responds by playing a (legal) vertex in $V(T_X)$ as her $\ell$th move. As before, Staller's move is an optimal move in the imaginary game played in $T$ by Lemma~\ref{lem1}. As his $\ell$th move in the real game, \D then plays an optimal move which he would play in the imaginary game in response to Staller's imagined $\ell$th move which plays a vertex in $V(T')$. From here onwards, all remaining moves of \D (in the imaginary game played in $T$) are played in the real game $T'$. Thus, once again, \D has a strategy to complete the game in $T'$ in at most~$\dstart(T) - 4$ moves. Therefore, $\dstart(T') \le \dstart(T) - 4$.

Every $\gt(T')$-set can be extended to a TD-set of $T$ by adding to it the four vertices $v$, $w$, $w'$ and $x$, implying that $\gt(T) \le \gt(T') + 4$. By assumption, $\gt(T) = \dstart(T)$. As observed earlier, $\dstart(T') \le \dstart(T) - 4$. Therefore,
\[
\dstart(T) = \gt(T) \le \gt(T') + 4 \le \dstart(T') + 4 \le \dstart(T).
\]
Consequently, we must have equality throughout the above inequality chain. In particular, $\gt(T') = \dstart(T')$. This completes the proof of Part~(b) of Lemma~\ref{lem7}.~\qed

\begin{lem}
\label{lem8}
If the leaf $u$ belongs to a $(3,3)$-structure, $T_x$, with pivot vertex~$x$, then some descendant of $x$ at distance~$3$ from $x$ belongs to a $(2,2)$-, $(3,1)$- or $(3,2)$-structure.
\end{lem}
\proof Let $T_x$ be a $(3,3)$-structure with pivot vertex~$x$ that contains the leaf $u$. Suppose, to the contrary, that no descendant of $x$ at distance~$3$ from $x$ belongs to a $(2,2)$-, $(3,1)$- or $(3,2)$-structure. Let $T_x$ be the path $uvwxw'v'u'$. Thus, $v$, $v'$, $w$ and $w'$ have degree~$2$ in $T$ and $u$ and $u'$ are leaves in $T$. By our earlier assumptions, $T$ has no strong support vertex. By supposition, no descendant of $x$ at distance~$3$ from $x$ belongs to a $(2,2)$-, $(3,1)$- or $(3,2)$-structure. This implies that the subtree of $T$ induced by $x$ and all its descendants can be obtained from a star with at least two leaves and with central vertex $x$ by subdividing every edge exactly twice.

If the root $r = x$, then the tree $T$ is obtained from a star on at least three vertices and with central vertex $x$ by subdividing every edge exactly twice. However, in this case, the vertex $x$ does not belong to any $\gamma_t(T)$-set, a contradiction. Therefore, the root $r \ne x$.

Suppose that the root $r = y$; that is, $r$ is the parent of $x$. By the Total Continuation Principle, the root $r$ has degree at least~$2$. Let $x'$ be a child of $r$ different from $x$. As her first move, \St plays the vertex $x'$. Since the vertex $x$ has at least two children, it is not possible for \D to block the vertex $x$ on his second move. Thus, on her second move, \St can play the vertex $x$. This, however, produces a contradiction since given the structure of the subtree of $T$ induced by $x$ and all its descendants, there is no $\gamma_t(T)$-set containing all three vertices $x$, $x'$ and $y$. Hence, the root $r \ne y$.

If the root $r$ is the grandparent of $x$, then \St plays the vertex $x$ as her first move. This produces a contradiction since there is no $\gamma_t(T)$-set containing both $x$ and its grandparent. If the root $r$ is not the grandparent of $x$, then \St plays the grandparent of $x$ as her first move. Since the vertex $x$ has at least two children, it is not possible for \D to block the vertex $x$ on his second move. Thus, on her second move, \St can play the vertex $x$. Once again, we produce a contradiction since there is no $\gamma_t(T)$-set containing both $x$ and its grandparent.

This final contradiction finishes the proof of the lemma since by choice of $x$ the root $r$ is either $x$ or some ancestor of $x$.~\qed

\begin{lem}
\label{lem9}
If the leaf $u$ belongs to a $(4)$-structure and $T' = T - \{u,v,w,x\}$, then $\gt(T') = \dstart(T')$.
\end{lem}
\proof Suppose that the leaf $u$ belongs to a $(4)$-structure. This implies that $v$, $w$ and $x$ all have degree~$2$ in $T$. We now consider the tree $T' = T - \{u,v,w,x\}$ and show that $\gt(T') = \dstart(T')$. As before, \D plays an imaginary game in $T$, although the real game is played in $T'$.
By the Total Continuation Principle, the first vertex played by \D in the imaginary game, namely $d_1 = r$, has degree at least~$2$ and is a legal move in the real game.

We show first that $\dstart(T') \le \dstart(T) - 2$. Suppose that every (optimal) move of \D in response to each of Staller's moves is a vertex in $T'$. In this case, every move of \St in the real game played in $T'$ is a legal move in the imaginary game played in $T$. We note that if \D played the vertex $y$ during the course of the game, then the vertex $y$ would have totally dominated at least one new vertex in the real game $T'$, for otherwise, by the Total Continuation Principle \D would have played the vertex $w$ instead. Thus, every move of \D is a legal move in the real game, implying that the game finishes in at most~$\dstart(T) - 2$ moves, since in the imaginary game played in $T$ at least two further moves are needed (in order to totally dominate the vertices $v$ and $w$). Thus, in this case, \D has a strategy to complete the game in $T'$ in at most~$\dstart(T) - 2$ moves. Hence, we may assume that during the imaginary game played in $T$, \D plays a vertex that does not belong to the real game $T'$.

By the Total Continuation Principle, the first vertex played by \D in the imaginary game that is not in $T'$ is not a leaf. Suppose that this occurs after Staller's $k$th move and that the vertex $z \in \{v,w\}$ is played by \D on his $(k+1)$st move. Let $\barz = \{v,w\} \setminus \{z\}$. \D imagines that \St immediately responds by playing $\barz$ as her $(k+1)$st move, which by Lemma~\ref{lem1}, is an optimal move in the imaginary game played in $T$. As his $(k+1)$st move in the real game, \D then plays an optimal move which he would play (as his $(k+2)$nd move) in the imaginary game in response to Staller's imagined move~$\barz$. By the Total Continuation Principle, we may assume that no subsequent move of \D plays the vertex $x$ which serves only to totally dominate the vertex $y$, since in this case \D would do at least as well by playing a neighbor of $y$ in $T'$. Thus all remaining moves of \D in the imaginary game played in $T$ are played in the real game $T'$. Further in this case, every move of \St in the real game played in $T'$ is a legal move in the imaginary game played in $T$. Thus, \D has a strategy to complete the game in $T'$ in at most~$\dstart(T) - 2$ moves.

Suppose therefore that \D plays as his $(k+1)$st move the vertex $x$ (and this is the first vertex played by \D in the imaginary game that is not in $T'$). \D imagines that \St immediately responds by playing $w$ as her $(k+1)$st move, which by Lemma~\ref{lem1}, is an optimal move in the imaginary game played in $T$.

Suppose that Dominator's $(k+2)$nd move in the imaginary game is the vertex $v$. Letting $S$ be the vertices totally dominated in $T$ after this move and applying Theorem~\ref{t:main1} to the partially total dominated forest $T|S$, we have
\[
\dstart(T) \stackrel{Lemma~\ref{lem1}}{=}  2k + 3 + \sstart(T|S) \stackrel{Theorem~\ref{t:main1}}{\ge}  2k + 3 + \dstart(T|S).
\]

\D now follows his optimal strategy in the partially total dominated forest $T|S$, and plays as his $(k+1)$st move in the real game his optimal first move in the game played in $T|S$. All subsequent moves of \D in response to Staller's moves are played in the real game $T'$ and are legal moves in the real game. Every subsequent move of \St in the real game played in $T'$ is a legal move in the imaginary game played in $T$, except possibly if she plays a vertex, $t$ say, that is a neighbor of $y$ and the only new vertex totally dominated by $t$ in the real game is the vertex $y$. Suppose that \St plays such a move $t$ as her $\ell$th move in the real game played in $T'$. \D now responds as follows. Immediately before she plays her move, we note that a total of $(2\ell -1) + 3$ moves are played in the imaginary game since, by our earlier assumptions, three additional moves $v$, $w$ and $x$ are played in the imaginary game. Letting $S'$ be the vertices totally dominated in $T$ after his $\ell$th move in the real game (namely, after his move that immediately precedes Staller's move $t$), and applying Theorem~\ref{t:main1} to the partially total dominated forest $T|S'$, we have
\[
\dstart(T) \stackrel{Lemma~\ref{lem1}}{=} 2\ell + 2  +  \sstart(T|S') \stackrel{Theorem~\ref{t:main1}}{\ge}  2\ell + 2 + \dstart(T|S').
\]

\D now follows his optimal strategy in the partially total dominated forest $T|S'$, and plays as his $(\ell+1)$st move in the real game his optimal first move in the game played in $T|S'$. All subsequent moves of \D and \St played in the real game $T'$ are legal moves. Thus, \D has a strategy to complete the game in $T'$ in at most~$\dstart(T) - 2$ moves, noting that three (redundant) moves were played in the imaginary game in $T$ (namely, the three vertices $v$, $w$ and $x$) and one additional move was played by \St (namely, the vertex $t$) that was not played in the imaginary game. Thus, once again \D has a strategy to complete the game in $T'$ in at most~$\dstart(T) - 2$ moves.

Hence, we may assume that Dominator's $(k+2)$nd move in the imaginary game belongs to $V(T')$ (and is therefore not the vertex $v$). In this case, he plays this move as his $(k+1)$st move in the real game. Continuing analogously as in the previous paragraphs, \D has a strategy to complete the game in $T'$ in at most~$\dstart(T) - 2$ moves. Thus, $\dstart(T') \le \dstart(T) - 2$.

Every $\gt(T')$-set can be extended to a TD-set of $T$ by adding to it the vertices $v$ and $w$, implying that $\gt(T) \le \gt(T') + 2$. By assumption, $\gt(T) = \dstart(T)$. As observed earlier, $\dstart(T') \le \dstart(T) - 2$. Therefore,
\[
\dstart(T) = \gt(T) \le \gt(T') + 2 \le \dstart(T') + 2 \le \dstart(T).
\]
Consequently, we must have equality throughout the above inequality chain. In particular, $\gt(T') = \dstart(T')$. This completes the proof of
Lemma~\ref{lem9}.~\qed

\section{The Family~$\cF$}
\label{S:cF}

In this section, we construct a family~$\cF^*$ of trees with equal total domination and game total domination number. For this purpose, we introduce some additional notation. Let $x$ be a specified vertex in a tree $T$. We define next several types of attachments at the vertex~$x$ that we use to build larger trees. In all cases, we call the vertex of the attachment that is joined to $x$ the \emph{link vertex} of the attachment. Recall that $P_n$ denotes a path on $n$ vertices.
\\[-20pt]
\begin{enumerate}
\item[$\bullet$] For $i \in [3]$, an \emph{attachment of Type-$i$ at $x$} is an operation that adds a path $P_{i+1}$ to $T$ and joins one of its ends to $x$.
\item[$\bullet$] An \emph{attachment of Type-$A$ at $x$} is obtained by adding an attachment of Type-$1$ at $x$ with link vertex $x'$, followed by at least one attachment of Type-$2$ at $x'$.
\item[$\bullet$] An \emph{attachment of Type-$B$ at $x$} is obtained by adding an attachment of Type-$A$ at $x$, followed by an attachment of Type-$3$ to at least one new (added) vertex at distance~$3$ from~$x$.
\end{enumerate}

We note that each attachment of Type-$A$ at $x$ can be obtained from a star $K_{1,k}$, for some $k \ge 2$, by subdividing $k-1$ edges twice and joining the central vertex of the original star to~$x$.

\subsection{The Family~$\cF_1$}
\label{S:cF1}

For integers $k_1,k_2,k_3,k_4 \ge 0$, let $\cT_{k_1,k_2,k_3,k_4}$ be the family of all trees  obtained from a trivial tree $K_1$ whose vertex is named~$a$ by applying $k_i$ attachments of Type-$i$ at $a$ for each $i \in [2]$, applying $k_3$ attachments of Type-$A$ at $a$ and applying $k_4$ attachments of Type-$B$ at $a$. Let
\[
\cF_1 = \bigcup_{k_1 \ge 1, k_2,k_3,k_4 \ge 0}
\cT_{k_1,k_2,k_3,k_4}.
\]

\begin{figure}[htb]
\begin{center}
\begin{tikzpicture}[scale=.7,style=thick,x=1cm,y=1cm]
\def\vr{2.5pt} 
\path (-0.5,10) coordinate (d1);
%
\path (-10,6) coordinate (a1);  \path (-8,6) coordinate (a2);
\path (-10,5) coordinate (a3);  \path (-8,5) coordinate (a4);
\path (-5,6) coordinate (b1); \path (-3,6) coordinate (b2);
\path (-5,5) coordinate (b3); \path (-3,5) coordinate (b4);
\path (-5,4) coordinate (b5); \path (-3,4) coordinate (b6);
\path (2,6) coordinate (c1); \path (4,6) coordinate (c2);
\path (1.5,5.2) coordinate (c3); \path (3.5,5.2) coordinate (c4);
\path (2,5) coordinate (c5); \path (4,5) coordinate (c6);
\path (2,4) coordinate (c7); \path (4,4) coordinate (c8);
\path (2,3) coordinate (c9); \path (4,3) coordinate (c10);
\path (7,6) coordinate (e1); \path (10,6) coordinate (e2);
\path (6.5,5.2) coordinate (e3); \path (9.5,5.2) coordinate (e4);
\path (7,5) coordinate (e5); \path (10,5) coordinate (e6);
\path (7,4) coordinate (e7); \path (10,4) coordinate (e8);
\path (6.5,3.2) coordinate (e9); \path (9.5,3.2) coordinate (e10);
\path (7,3) coordinate (e11); \path (10,3) coordinate (e12);
\path (7,2) coordinate (e13); \path (10,2) coordinate (e14);
\path (7,1) coordinate (e15); \path (10,1) coordinate (e16);
\path (6.5,.2) coordinate (e17); \path (9.5,.2) coordinate (e18);

%
\draw (d1) -- (a1); \draw (d1) -- (a2); \draw (d1) -- (b1); \draw (d1) -- (b2);
\draw (d1) -- (c1); \draw (d1) -- (c2); \draw (d1) -- (e1); \draw (d1) -- (e2);
\draw (a1) -- (a3); \draw (a2) -- (a4);
\draw (b1) -- (b3); \draw (b2) -- (b4); \draw (b3) -- (b5); \draw (b4) -- (b6);
\draw (c1) -- (c3); \draw (c2) -- (c4); \draw (c1) -- (c5); \draw (c2) -- (c6);
\draw (c5) -- (c7); \draw (c6) -- (c8); \draw (c7) -- (c9); \draw (c8) -- (c10);
\draw (e1) -- (e3); \draw (e2) -- (e4); \draw (e1) -- (e5); \draw (e2) -- (e6);
\draw (e5) -- (e7); \draw (e6) -- (e8); \draw (e7) -- (e9); \draw (e8) -- (e10);
\draw (e7) -- (e11); \draw (e8) -- (e12); \draw (e11) -- (e13); \draw (e12) -- (e14);
\draw (e13) -- (e15); \draw (e14) -- (e16); \draw (e15) -- (e17); \draw (e16) -- (e18);

\draw [fill=white] (-.65,9.85) rectangle (-.35,10.15);
\draw (a1) [fill=white] circle (\vr); \draw (a2) [fill=black] circle (\vr);
\draw (a3) [fill=white] circle (\vr); \draw (a4) [fill=white] circle (\vr);
\draw (b1) [fill=white] circle (\vr); \draw (b2) [fill=white] circle (\vr);
\draw (b3) [fill=white] circle (\vr); \draw (b4) [fill=white] circle (\vr);
\draw (b4) [fill=white] circle (\vr);
\draw (b5) [fill=white] circle (\vr); \draw (b6) [fill=white] circle (\vr);
\draw (c1) [fill=black] circle (\vr); \draw (c2) [fill=black] circle (\vr);
\draw (c3) [fill=white] circle (\vr); \draw (c4) [fill=white] circle (\vr);
\draw (c5) [fill=white] circle (\vr); \draw (c6) [fill=white] circle (\vr);
\draw (c7) [fill=white] circle (\vr); \draw (c8) [fill=white] circle (\vr);
\draw (c9) [fill=white] circle (\vr); \draw (c10) [fill=white] circle (\vr);
\draw (e1) [fill=black] circle (\vr); \draw (e2) [fill=black] circle (\vr);
\draw (e3) [fill=white] circle (\vr); \draw (e4) [fill=white] circle (\vr);
\draw (e5) [fill=white] circle (\vr); \draw (e6) [fill=white] circle (\vr);
\draw (e7) [fill=white] circle (\vr); \draw (e8) [fill=white] circle (\vr);
\draw (e9) [fill=white] circle (\vr); \draw (e10) [fill=white] circle (\vr);
\draw (e11) [fill=white] circle (\vr); \draw (e12) [fill=white] circle (\vr);
\draw (e13) [fill=white] circle (\vr); \draw (e14) [fill=white] circle (\vr);
\draw (e15) [fill=white] circle (\vr); \draw (e16) [fill=white] circle (\vr);
\draw (e17) [fill=white] circle (\vr); \draw (e18) [fill=white] circle (\vr);
\draw (-0.5,10.5) node {$a$};
\draw (-9,4.5) node {$\underbrace{\phantom{1111111}}$};
\draw (-9,3.8) node {$k_1\ge 1$};
\draw (-4,3.5) node {$\underbrace{\phantom{1111111}}$};
\draw (-4,2.8) node {$k_2 \ge 0$};
\draw (3,2.5) node {$\underbrace{\phantom{1111111}}$};
\draw (3,1.8) node {$k_3 \ge 0$};
\draw (2,4) node [draw,scale=0.55,diamond,fill=white!50]{};
\draw (4,4) node [draw,scale=0.55,diamond,fill=white!50]{};
\draw (8,-.3) node {$\underbrace{\phantom{1111111111}}$};
\draw (8,-1.0) node {$k_4 \ge 0$};
\end{tikzpicture}
\end{center}
\vskip -0.75 cm
\caption{A tree in the family~$\cT_{k_1,k_2,k_3,k_4}$} \label{fig:F1}
\end{figure}
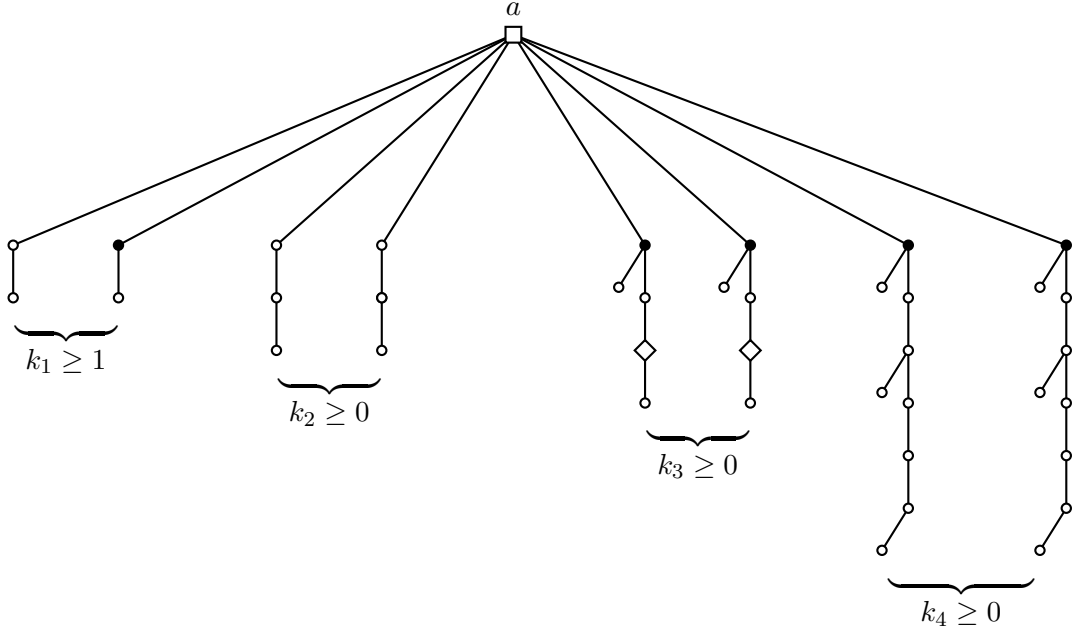

A tree $T$ in the family~$\cT_{k_1,k_2,k_3,k_4}$ is illustrated in Figure~\ref{fig:F1}. We note that there can be additional attachments of type-$2$ at each darkened vertex in Figure~\ref{fig:F1} that belongs to an attachment of type-$A$. We shall show (see the proof of Lemma~\ref{lem10}) that the vertex $a$ (depicted by the open square in Figure~\ref{fig:F1}) is an optimal first move of Dominator.

\begin{lem}
\label{lem10}
If $T \in \cF_1$, then $\gt(T) = \dstart(T)$.
\end{lem}
\proof Let $T \in \cF_1$, and so $T = T_{k_1,k_2,k_3,k_4}$, for some integers $k_1 \ge 1$ and $k_2,k_3,k_4 \ge 0$. Adopting our earlier notation, let $a$ be the vertex of the trivial tree $K_1$ from which the tree $T$ was built. We first consider the case when every link vertex of an attachment of Type-$A$ or Type-$B$ has degree exactly~$3$. Thus, every attachment of Type-$A$ at $a$ is obtained by adding an attachment of Type-$1$ at $a$ with link vertex $x'$, followed by exactly one attachment of Type-$2$ at $x'$.

Every TD-set of $T$ necessarily contains all its support vertices and a neighbor of each support vertex, implying that such a set contains at least one vertex from every Type-$1$ attachment, two vertices from every Type-$2$ attachment, three vertices from every Type-$A$ attachment, and five vertices from every Type-$B$ attachment. Further, a TD-set of $T$ that contains only one vertex from some Type-$1$ attachment, also contains the vertex~$a$. It follows that $\gt(T) \ge k_1 + 2k_2 + 3k_3 + 5k_4 + 1$.  Conversely, the set consisting of all vertices of $T$ that are not leaves, and that are not degree-$2$ vertices at distance~$4$ from $a$ in $T$, forms a TD-set of $T$ of size~$k_1 + 2k_2 + 3k_3 + 5k_4 + 1$, and so $\gt(T) \le k_1 + 2k_2 + 3k_3 + 5k_4 + 1$. Consequently, $\gt(T) = k_1 + 2k_2 + 3k_3 + 5k_4 + 1$.

We show next that \D has a strategy to finish the game in $k_1 + 2k_2 + 3k_3 + 5k_4 + 1$ moves. \D plays the vertex $a$ as his first move. This first move of \D  blocks \St from playing a leaf at distance~$2$ from $a$, implying that exactly one vertex is played from every attachment of Type-$1$. Further, this first move of \D implies that exactly three vertices are played in every attachment of Type-$A$. \D now adopts the following strategy.

If \St plays a leaf (at distance~$3$ from $a$) in an attachment of Type-$2$, then \D responds as follows. If no support vertex in an attachment of Type-$1$ has yet been played, then \D plays such a neighbor of $a$. Otherwise, \D plays any playable vertex that is not a leaf, playing a support vertex wherever possible. This strategy of \D implies that exactly two vertices are played in every attachment of Type-$2$.

Finally, suppose that $F$ is an attachment of Type-$B$, where $F$ is the path $v_1v_2 \ldots v_8$ with a pendant edge $v_4u_4$ and with link vertex $v_2$ (and so, $av_2$ is an edge of $T$). As observed earlier, the leaf $v_1$ is not playable. If \St plays one of the neighbors of $v_4$, namely one of the vertices $v_3$, $v_5$ or $u_4$, then \D immediately responds by playing the vertex $v_7$ if it has not yet been played; otherwise he plays any playable vertex that is not a leaf, playing a support vertex whenever possible. If \St plays one of the neighbors of $v_7$, namely one of the vertices $v_6$ or $v_8$, then \D immediately responds by playing the vertex $v_4$ if it has not yet been played; otherwise he plays any playable vertex that is not a leaf, playing a support vertex whenever possible. If \St plays some other vertex from $F$, then \D plays any playable vertex that is not a leaf, playing a support vertex wherever possible. This strategy of \D implies that exactly five vertices are played in every attachment of Type-$B$. Thus, \D has a strategy to finish the game in $k_1 + 2k_2 + 3k_3 + 5k_4 + 1$ moves. Hence, $\dstart(T) \le k_1 + 2k_2 + 3k_3 + 5k_4 + 1 = \gt(T)$. Since $\gt(T) \le \dstart(T)$, this implies that $\gt(T) = \dstart(T)$.
An analogous proof works if we relax the requirement that every link vertex of an attachment of Type-$A$ or Type-$B$ has degree exactly~$3$.~\qed

\subsection{The Tree $F_{10}$}

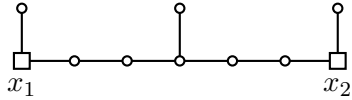
\begin{figure}[htb]
\begin{center}
\begin{tikzpicture}[scale=.7,style=thick,x=1cm,y=1cm]
\def\vr{2.5pt} 
\path (-3,1) coordinate (a1);  \path (-3,0) coordinate (a2);
\path (-2,0) coordinate (a3);  \path (-1,0) coordinate (a4);
\path (0,0) coordinate (a5);  \path (1,0) coordinate (a6);
\path (2,0) coordinate (a7);  \path (3,0) coordinate (a8);
\path (3,1) coordinate (a9);  \path (0,1) coordinate (a10);
\draw (a1) -- (a2); \draw (a2) -- (a3);
\draw (a3) -- (a4); \draw (a4) -- (a5);
\draw (a5) -- (a6); \draw (a6) -- (a7);
\draw (a7) -- (a8); \draw (a8) -- (a9);
\draw (a5) -- (a10);

\draw (a1) [fill=white] circle (\vr); \draw [fill=white] (-3.15,-.15) rectangle (-2.85,.15);
\draw (a3) [fill=white] circle (\vr); \draw (a4) [fill=white] circle (\vr);
\draw (a5) [fill=white] circle (\vr); \draw (a6) [fill=white] circle (\vr);
\draw (a7) [fill=white] circle (\vr); \draw [fill=white] (2.85,-.15) rectangle (3.15,.15);
\draw (a9) [fill=white] circle (\vr); \draw (a10) [fill=white] circle (\vr);
\draw (-3,-.5) node {$x_1$};
\draw (3,-.5) node {$x_2$};
\end{tikzpicture}
\end{center}
\vskip -0.5 cm
\caption{The tree $F_{10}$} \label{fig:F3}
\end{figure}

Let $F_{10}$ be the tree of order~$10$ obtained from a star $K_{1,3}$ by subdividing two edges three times. The tree $F_{10}$ is illustrated in Figure~\ref{fig:F3}. The two vertices $x_1$ and $x_2$ (represented by an open square) are the two optimal first moves of Dominator.

\begin{lem}
\label{lem12}
If $T = F_{10}$, then $\gt(T) = \dstart(T)$.
\end{lem}
\proof Let $F_{10}$ be the path $v_1v_2 \ldots v_9$ with a pendant edge $v_5u_5$. Every TD-set of $T$ necessarily contains all its support vertices and a neighbor of each support vertex, implying that $\gt(T) = 6$. We show next that \D has a strategy to finish the game in six moves. \D plays the vertex $v_2$ as his first move, and adopts the following strategy. If \St plays one of the neighbors of $v_5$, then \D immediately responds by playing the vertex $v_8$, if it has not yet been played. If \St plays one of the neighbors of $v_2$ or $v_8$, then \D immediately responds by playing the vertex $v_5$, if it has not yet been played. Otherwise, \D plays any playable vertex that is not a leaf, playing a support vertex wherever possible. This strategy of \D implies that exactly six vertices are played. Thus, \D has a strategy to finish the game in six moves. Hence, $\dstart(T) \le 6 = \gt(T)$. Since $\gt(T) \le \dstart(T)$, this implies that $\gt(T) = \dstart(T)$.~\qed

\subsection{The Family $\cF^*$}

Let $\cF = \cF_1 \cup \{K_2,F_{10}\}$. By Lemmas~\ref{lem10} and~\ref{lem12}, and since %
$\gt(K_2) = \dstart(K_2) = 2$,
we have the following result.

\begin{lem}
\label{lem13}
If $T \in \cF$, then $\gt(T) = \dstart(T)$.
\end{lem}

Let $\cF^*$ be the family of all stars on at least two vertices together with all trees that can be obtained from a tree $F$ of order at least~$3$ in the family~$\cF$ by adding any number, including the possibility of zero, additional pendant edges to support vertices of $F$. As a consequence of Lemmas~\ref{lem_strong} and~\ref{lem13}, and the observation that $\gt(T) = \dstart(T) = 2$ for every star on at least two vertices, we have the following result.

\begin{lem}
\label{lem14}
If $T \in \cF^*$, then $\gt(T) = \dstart(T)$.
\end{lem}

\section{Proof of Theorem~\ref{t:main2}}
\label{S:proof2}

In this section we present a proof of Theorem~\ref{t:main2}. Recall its statement.

\noindent \textbf{Theorem~\ref{t:main2}}. \emph{Let $T$ be a nontrivial tree. Then, $\gt(T) = \dstart(T)$ if and only if $T \in \cF^*$.}

\medskip
\proof The sufficiency follows from Lemma~\ref{lem14}. To prove the necessity, we proceed by induction on the order~$n \ge 2$ of a tree $T$ that satisfies $\gt(T) = \dstart(T)$. If $n \in \{2,3\}$, then $T$ is a star, and so $T \in \cF^*$. This establishes the base case. Let $n \ge 4$ and suppose that every nontrivial tree $T'$ of order~$n'$, where $2 \le n' < n$, satisfying $\gt(T') = \dstart(T')$ belongs to the family~$\cF^*$. Let $T$ be a tree of order~$n$ that satisfies $\gt(T) = \dstart(T)$. If $T$ is a star, then $T \in  \cF^*$. Hence, we may assume that $\diam(T) \ge 3$. If $\diam(T) = 3$, then $\gt(T) = 2$ and $\dstart(T) = 3$, a contradiction. Hence, $\diam(T) \ge 4$. By Lemma~\ref{lem2}, no two leaves are at distance~$3$ apart in $T$.

We now root the tree $T$ at an optimal first move, $r$ say, of Dominator. Thus, \D plays the vertex $d_1 = r$ as his first move in the dominator-start total game played in $T$. By the Total Continuation Principle, it is never in Dominator's best interests to play a leaf. Hence, $d_T(r) \ge 2$. Let $u$ be a vertex at maximum distance from the root $r$ in $T$. Necessarily, $u$ is a leaf. Let $v$ be the parent of $u$ and let $w$ the parent of $v$. Further, if $r \ne w$, let $x$ the parent of $w$, and if $r \ne x$, let $y$ the parent of $x$.

\begin{unnumbered}{Claim~\ref{t:main2}.1}
If $T$ contains a strong support vertex, then $T \in \cF^*$. \end{unnumbered}
\ \\[-2em]
\noindent\textbf{Proof.}
Suppose that $T$ contains a strong support vertex. Thus, $T$ contains a $(1,1)$-structure, $T_v$, with pivot vertex~$v$. Let $T'$ be the tree obtained from $T$
by deleting a leaf vertex in $T_v$. By Lemma~\ref{lem3}, $\gt(T') = \dstart(T')$. Applying the inductive hypothesis to $T'$, the tree $T' \in \cF^*$. By definition of the family $\cF^*$, every tree that can be obtained from a tree in the family~$\cF^*$ by adding an additional pendant edge to a support vertex also belongs to $\cF^*$. In particular, the tree $T \in \cF^*$.~\smallqed

\medskip
By Claim~\ref{t:main2}.1, we may assume that $T$ contains no strong support vertex, for otherwise the desired result follows. In particular, $d_T(v) = 2$. More generally, the parent of a vertex at maximum distance from the root $r$ in $T$ has degree~$2$.

\begin{unnumbered}{Claim~\ref{t:main2}.2}
If $d_T(w) \ge 3$, then $T \in \cF$.
\end{unnumbered}
\ \\[-2em]
\noindent\textbf{Proof.} Suppose that $d_T(w) \ge 3$. By Lemma~\ref{lem4}, the tree $T$ contains no $(2,1)$-structure. Hence, no child of $w$ is a leaf. Since $d_T(w) \ge 3$, this implies that the leaf $u$ belongs to a $(2,2)$-structure. By Lemma~\ref{lem5}, the tree $T$ is a subdivided star. Thus, $T = T_{k_1,0,0,0}$ for some integer $k_1 \ge 2$, and so $T \in \cF_1 \subset \cF$.~\smallqed

\medskip
By Claim~\ref{t:main2}.2, we may assume that $d_T(w) = 2$, for otherwise the desired result follows. More generally, we may assume that the grandparent of a vertex at maximum distance from the root $r$ in $T$ has degree~$2$.

\begin{unnumbered}{Claim~\ref{t:main2}.3}
If $d_T(x) \ge 3$ and a child of $x$ is a leaf, then $T \in \cF$.
\end{unnumbered}
\ \\[-2em]
\noindent\textbf{Proof.} Suppose that $d_T(x) \ge 3$ and that a child of $x$ is a leaf. In this case, the leaf $u$ belongs to a $(3,1)$-structure. Let $T' = T - \{u,v,w\}$. We note that $x$ is a support vertex of $T'$. By Lemma~\ref{lem6}, $\gt(T') = \dstart(T')$. Further, $r = y$, where we recall that $y$ is the parent of $x$. This implies that $r$ is within distance~$4$ from every vertex of $T$ and therefore $\diam(T') \le \diam(T) \le 8$. As observed earlier, no two leaves are at distance~$3$ apart in $T$. Thus, since $d_T(r) \ge 2$ and $x$ is a support vertex, we note that $\diam(T) \ge \diam(T') \ge 4$. In particular, $n \ge 8$. Applying the inductive hypothesis to $T'$, the tree $T' \in \cF^*$.  Since $T$ has no strong support vertex, neither does the tree $T'$. Thus, $T' \in \cF$.

Suppose that $T' = F_{10}$. If the support vertex $x$ is the vertex named $x_1$ (or by symmetry the vertex named $x_2$) in Figure~\ref{fig:F3}, then $\diam(T) = 10$, a contradiction. If the support vertex $x$ is the central vertex of the tree $T' = F_{10}$, then there is a vertex in $T'$ at distance~$5$ from the vertex $r = y$, a contradiction. Since both cases produce a contradiction, $T' \ne F_{10}$.

Hence, $T' \in \cF_1$. Thus, $T' \in \cT_{k_1,k_2,k_3,k_4}$ for some integers $k_1 \ge 1$ and $k_2,k_3,k_4 \ge 0$. Let $r'$ be the vertex in the trivial tree $K_1$ used to build the tree $T'$ (and so, $r'$ corresponds to the vertex named ``$a$" in Figure~\ref{fig:F1}). If $k_4 > 0$, then $\diam(T') > 8$, and so $\diam(T) > 8$, a contradiction. Hence, $k_4 = 0$ and $T' \in \cT_{k_1,k_2,k_3,0}$.

Suppose that $x$ is the link vertex of a Type-$1$ attachment in $T'$. Suppose that $k_1 = 1$. In this case, $T \in \cT_{0,k_2,k_3+1,0}$. The structure of the tree $T$ implies that the vertex $r'$ belongs to no $\gt(T)$-set. However, $n \ge 8$, and so $k_2 + k_3 + 1 \ge 2$, implying that on Dominator's first move, he cannot block \St from playing the vertex~$r'$ as her first move. Thus, \St has a strategy to force at least~$\gt(T) + 1$ moves in the game, a contradiction. Hence, $k_1 \ge 2$, implying that $T \in \cT_{k_1 - 1,k_2,k_3+1,0} \subset \cF_1$.

Suppose that $x$ belongs to a Type-$2$ attachment in $T'$ (and is therefore a support vertex at distance~$2$ from $r'$ in $T'$). Let $w_1w_2w_3$ be the Type-$2$ attachment at $r'$ that contains $x$, where $w_1$ is the link vertex (and so, $w_1$ is joined to $r'$) and $x = w_2$. Let $w_4w_5w_6$ be the attachment of Type-$2$ at $x$ with link vertex $w_4$, and so $w_4$ is joined to $w_2$. Let $r'y_1y_2$ be a path emanating from $r'$, where the path $y_1y_2$ represents a Type-$1$ attachment in $T'$ at $r'$ with link vertex $y_1$. We note that $y_2y_1r'w_1w_2w_4w_5w_6$ is a path in $T$ and that $w_2w_3$ is a pendant edge in $T$. Further, each of $y_1,w_1,w_4,w_5$ have degree~$2$ in $T$ and each of $y_2, w_3,w_6$ are leaves in $T$, while $d_T(r') \ge 2$. The structure of $T$ implies that neither vertex $w_1$ nor $w_3$ belongs to a $\gt(T)$-set. However, on Dominator's first move he cannot block \St from playing one of the vertices $w_1$ nor $w_3$ as her first move. Thus, \St has a strategy to force at least~$\gt(T) + 1$ moves in the game, a contradiction. Hence, $x$ does not belong to a Type-$2$ attachment in $T'$.

Suppose, finally, that the support vertex $x$ belongs to a Type-$A$ attachment in $T'$. Thus, $x$ is either the link vertex of the attachment (that is adjacent to $r'$) or the vertex in the attachment at distance~$3$ from $r'$ in $T'$. If $x$ is the link vertex of the attachment, then $T \in \cT_{k_1,k_2,k_3,0} \subset \cF_1$. If $x$ is the vertex in the attachment at distance~$3$ from $r'$ in $T'$, then $T \in \cT_{k_1,k_2,k_3,1} \subset \cF_1$. In both cases, $T \in \cF_1$.~\smallqed

\begin{unnumbered}{Claim~\ref{t:main2}.4}
If $d_T(x) \ge 3$ and a grandchild of $x$ is a leaf, then $T \in \cF$.
\end{unnumbered}
\ \\[-2em]
\noindent\textbf{Proof.} Suppose that $d_T(x) \ge 3$ and a grandchild, $v'$, of $x$ is a leaf. Recall that $u$ is a vertex at maximum distance from the root $r$ in $T$. Further, recall that $T$ contains no strong support vertex and no two leaves are at distance~$3$ apart in $T$. These observations imply that the parent, $w'$ say, of the grandchild $v'$ has degree~$2$ in $T$. Thus, the leaf $u$ belongs to a $(3,2)$-structure, $T_x$, with pivot vertex~$x$. We note that $T_x$ is the path $uvwxw'v'$, where $v,w,w'$ have degree~$2$ in $T$ and $u$ and $v'$ are leaves in $T$.

If the root $r = x$, then, by Lemma~\ref{lem7}(a), $T$ is obtained from a star with $k_1 + k_2$ leaves, where $k_1,k_2 \ge 1$, by subdividing $k_1$ edges once and $k_2$ edges twice. Thus, $T \in \cT_{k_1,k_2,0,0} \subset \cF_1$. Hence, we may assume that the root $r \ne x$, for otherwise the desired result follows. With this assumption, Lemma~\ref{lem7}(b) implies that the vertex $x$ has degree~$3$. Further, the root $r$ is the grandparent of $x$ (equivalently, the root $r$ is the parent of $y$). The $(3,2)$-structure $T_x$, is illustrated in Figure~\ref{f:S32}. We note that the parent $y$ of $x$ is within distance~$6$ from every vertex of $T$, and the root $r$ is within distance~$5$ from every vertex of $T$. In particular, $\diam(T') \le \diam(T) \le 10$. We now consider the tree $T' = T - V(T_x)$. By Lemma~\ref{lem7}(b), $\gt(T') = \dstart(T')$. Applying the inductive hypothesis to $T'$, the tree $T' \in \cF^*$.

By assumption, $T$ has no strong support vertex. Suppose that $T'$ contains a strong support vertex. Necessarily, such a strong support vertex is the root $r$. Let $y'$ be the leaf-neighbor of $r$ different from $y$ in $T'$. If $d_T(r) = 2$, then the tree $T$ is determined and $T \in \cT_{1,2,0,0} \subset \cF_1$. Hence, we may assume that $d_T(r) \ge 3$. Let $y''$ be a neighbor of $r$ different from $y$ and $y'$. Since $T$ has no strong support vertex, the vertex $y''$ has degree at least~$2$. Recall that \D plays the vertex $r$ as his first move in the dominator-start total game played in $T$. \St responds by playing the leaf $y'$ as her first move. \D is now unable to block both vertices $y$ and $y''$, and on Staller's second move she plays one of these two vertices. Thus, \St has a strategy to force at least~$\gt(T) + 1$ moves in the game, a contradiction. Hence, $T'$ contains no strong support vertex. Thus, $T' \in \cF$, and so $T' \in \{F_{10}\} \cup \cF_1$.

Suppose that $T' = F_{10}$. Let $c$ denote the central vertex of $T'$. As observed earlier, the vertices $r$ and $y$ are within distance~$5$ and~$6$, respectively, from every vertex of $T$. This implies that $r = c$ or $r$ is a neighbor of $c$ in $T'$. Further, the vertex $y$ is within distance~$2$ from $c$ in $T'$. A simple case analysis, noting that \D plays the vertex $r$ as his first move in $T$ shows that \St has a strategy to force at least three vertices played in the closed neighborhood of $c$ in $T$, thereby forcing at least~$\gt(T) + 1$ moves in the game, a contradiction.
If $T' \in \cF_1$, then a tedious, but straightforward, analysis shows that \St has a strategy to force at least~$\gt(T) + 1$ moves in the game, a contradiction.~\smallqed

\medskip
By Claim~\ref{t:main2}.3 and Claim~\ref{t:main2}.4, we may assume that if $d_T(x) \ge 3$, then no child and no grandchild of $x$ is a leaf, for otherwise $T \in \cF$, as desired. With this assumption, we have the following claim.

\begin{unnumbered}{Claim~\ref{t:main2}.5}
$d_T(x) = 2$.
\end{unnumbered}
\ \\[-2em]
\noindent\textbf{Proof.} Suppose that $d_T(x) \ge 3$. By assumption, no child and no grandchild of $x$ is a leaf. Thus, no descendant of $x$ at distance~$3$ from $x$ belongs to a $(2,2)$-, $(3,1)$- or $(3,2)$-structure. This contradicts Lemma~\ref{lem8}.~\smallqed

\medskip
By Claim~\ref{t:main2}.5, $d_T(x) = 2$. Recall that by our earlier assumptions, $d_T(v) = d_T(w) = 2$. Thus, the leaf $u$ belongs to a $(4)$-structure. Let $T' = T - \{u,v,w,x\}$. By Lemma~\ref{lem9}, $\gt(T') = \dstart(T')$. Applying the inductive hypothesis to $T'$, the tree $T' \in \cF^*$.
By assumption, $T$ has no strong support vertex. Suppose that $T'$ contains a strong support vertex. Necessarily, such a strong support vertex is the parent, $z$ say, of $y$. Thus, $d_T(y) = 2$ and $z$ is a support vertex in $T$. \St plays as her first move in $T$ the vertex~$x$, thereby forcing at least~$\gt(T) + 1$ moves in the game, a contradiction. Hence, $T'$ contains no strong support vertex. Thus, $T' \in \cF$, and so $T' \in \{F_{10}\} \cup \cF_1$. If $T' = F_{10}$, then a tedious, but straightforward, analysis shows that \St has a strategy to force at least~$\gt(T) + 1$ moves in the game, a contradiction.

Hence, $T' \in \cF_1$. Thus, $T' \in \cT_{k_1,k_2,k_3,k_4}$ for some integers $k_1 \ge 1$ and $k_2,k_3,k_4 \ge 0$. Let $r'$ be the vertex in the trivial tree $K_1$ used to build the tree $T'$ (and so, $r'$ corresponds to the vertex named ``$a$" in Figure~\ref{fig:F1}). A tedious, but straightforward, analysis shows that if the vertex $y$ is not a support vertex of degree~$2$ in $T'$ at distance~$3$ from $r'$ (such a vertex is depicted by a diamond in Figure~\ref{fig:F1}), then  \St has a strategy to force at least~$\gt(T) + 1$ moves in the game, a contradiction. Thus, the vertex $y$ is a support vertex of degree~$2$ in $T'$ at distance~$3$ from $r'$. Thus, $T \in \cT_{k_1,k_2,k_3-1,k_4 + 1} \subset \cF_1$. This completes the proof of Theorem~\ref{t:main2}.~\qed

\section*{Acknowledgements}

The first author is supported in part by the South African National Research Foundation and the University of Johannesburg. Research of the second author is supported by a grant from the Simons Foundation (\#209654 to Douglas Rall).

\end{document}